\documentclass[11pt,a4paper,fullpage]{article}
\usepackage[english]{babel}
\usepackage[utf8]{inputenc}
\usepackage{fancyhdr}
\usepackage[numbers]{natbib}

\usepackage{graphicx}
\usepackage{amsmath}
\usepackage{extsizes}
\usepackage{relsize}
\usepackage{multicol}
\usepackage{ragged2e}
\usepackage{fontenc}
\usepackage{mathcomp}
\usepackage{latexsym}
\usepackage{amssymb}
\usepackage{times}
\newtheorem{Teorema}{Theorem}[section]

\newtheorem{Posledica}[Teorema]{Corollary}

\newtheorem{Lema}[Teorema]{Lemma}
\newtheorem{Primedba}[Teorema]{Remark}
\usepackage{dsfont}
\newcommand\scalemath[2]{\scalebox{#1}{\mbox{\ensuremath{\displaystyle #2}}}} 
\numberwithin{equation}{section}
\hyphenchar\font=-1
\begin{document}
	\title {Invertibility properties of operator matrices on Hilbert spaces}	

\author{Nikola Sarajlija\footnote{corresponding author: Nikola Sarajlija; University of Novi Sad, Faculty of Sciences, Novi Sad 21000, Serbia; {\it e-mail}: {\tt nikola.sarajlija@dmi.uns.ac.rs}}\footnote{The author is supported by the Ministry of Education, Science and Technological Development of the Republic of Serbia under grant no.  451-03-68/2022-14/200125.}}
\maketitle

\begin{abstract}
Denote by $T_n^d(A)$ an upper triangular operator matrix of dimension $n$ whose diagonal entries are given and the others are unknown. In this article we provide necessary and sufficient conditions for various types of Fredholm and Weyl completions of $T_n^d(A)$.  As consequences, we get corrected perturbation results of Wu et al. (Ann. Funct. Anal. \textbf{11} (2020), no. 3, 780–798 and Acta Math. Sin. (Engl. Ser.) \textbf{36} (2020), no. 7, 783–796). In the special case $n=2$, we recover many already existing known results, and specially we correct results of Zhang et al. (J. Math. Anal. Appl. \textbf{392} (2012), no. 2, 103–110). Finally, in the case of essential Fredholm invertibility, in the special case $n=2$ we obtain some results that seem new in the literature. Our method is based on the space decomposition technique, similarly to the work of Huang et al. (Math. Nachr. \textbf{292} (2019), no. 11, 2411–2426 and Math. Nachr. \textbf{291} (2018), no. 1, 187–203), but our approach extends to arbitrary dimension $n\geq2$. 
\end{abstract}

\textit{$2020$ Math. Subj. Class:} 47A08, 47A53, 47A55, 47A05, 47A10.

\vspace{2mm}

\section{Preliminaries}
Throughout this text we always assume that $\mathcal{H},\mathcal{K},\mathcal{H}_1,...,\mathcal{H}_n$ are infinite dimensional separable Hilbert spaces. Collection of operators (linear and bounded mappings) from $\mathcal{H}$ to $\mathcal{K}$ is denoted by $\mathcal{B}(\mathcal{H},\mathcal{K})$, where $\mathcal{B}(\mathcal{H}):=\mathcal{B}(\mathcal{H},\mathcal{H})$. Letter $T$ is most commonly used to denote some operator, and we keep standard notation $\mathcal{N}(T)$, $\mathcal{R}(T)$ for the null and range space of $T$, respectively. As always, if $T\in\mathcal{B}(\mathcal{H},\mathcal{K})$, then $T^*\in\mathcal{B}(\mathcal{K},\mathcal{H})$ is the adjoint of $T$.

In this article we investigate various types of invertibility that occur in Fredholm theory. First, we introduce some notation \cite{ZANA}. Let $T\in\mathcal{B}(\mathcal{H})$, $\alpha(T)=\dim\mathcal{N}(T)$ and $\beta(T)=\dim\mathcal{N}(T^*)$. In the previous sentence, $\dim$ is the orthogonal Hilbert dimension. We call $\alpha(T)$ and $\beta(T)$ the nullity and the deficiency of $T$, respectively, and in the case where at least one of them is finite we define $\mathrm{ind}(T)=\alpha(T)-\beta(T)$ to be the index of $T$. Now we are able to define families of upper and lower semi-Fredholm operators, respectively, as
$$
\begin{aligned}
\Phi_+(\mathcal{H})=\left\{ T\in\mathcal{B}(\mathcal{H}):\ \alpha(T)<\infty\ and\ \mathcal{R}(T)\ is \ closed\right\}
\end{aligned}
$$
and
$$
\begin{aligned}
\Phi_-(\mathcal{H})=\lbrace T\in\mathcal{B}(\mathcal{H}): \beta(T)<\infty\rbrace.
\end{aligned}
$$
The set of Fredholm operators is
$$\Phi(\mathcal{H})=\Phi_+(\mathcal{H})\cap\Phi_-(\mathcal{H})=\lbrace T\in\mathcal{B}(\mathcal{H}): \alpha(T)<\infty\ and\ \beta(T)<\infty\rbrace.$$
Put
$$\Phi_+^-(\mathcal{H})=\lbrace T\in\Phi_+(\mathcal{H}):\ \mathrm{ind}(T)\leq0\rbrace$$
and
$$\Phi_-^+(\mathcal{H})=\lbrace T\in\Phi_-(\mathcal{H}):\ \mathrm{ind}(T)\geq0\rbrace.$$
These are collections of upper and lower semi-Weyl operators, respectively.

Corresponding spectra of an operator $T\in\mathcal{B}(\mathcal{H})$ are defined as follows:\\
-the upper semi-Fredholm spectrum: $\sigma_{SF+}(T)=\lbrace\lambda\in\mathds{C}: \lambda-T\not\in\Phi_{+}(\mathcal{H})\rbrace$;\\
-the lower semi-Fredholm spectrum: $\sigma_{SF-}(T)=\lbrace\lambda\in\mathds{C}: \lambda-T\not\in\Phi_{-}(\mathcal{H})\rbrace$;\\
-the essential spectrum: $\sigma_{e}(T)=\lbrace\lambda\in\mathds{C}: \lambda-T\not\in\Phi(\mathcal{H})\rbrace$;\\
-the upper semi-Weyl spectrum: $\sigma_{aw}(T)=\lbrace\lambda\in\mathds{C}: \lambda-T\not\in\Phi_+^-(\mathcal{H})\rbrace$;\\
-the lower semi-Weyl spectrum: $\sigma_{sw}(T)=\lbrace\lambda\in\mathds{C}: \lambda-T\not\in\Phi_-^+(\mathcal{H})\rbrace$;\\[1mm]
All of these spectra are compact nonempty subsets of the complex plane.\\
We write $\rho_{SF+}(T), \rho_{SF-}(T), \rho_e(T), \rho_{aw}(T), \rho_{sw}(T)$ to denote their complements, respectively.

Let $D_1\in\mathcal{B}(\mathcal{H}_1),\ D_2\in\mathcal{B}(\mathcal{H}_2),...,D_n\in\mathcal{B}(\mathcal{H}_n)$ be given. Partial upper triangular operator matrix of dimension $n$ is \cite{WU2}, \cite{WU3}
\begin{equation}\label{OSNOVNI}
T_n^d(A)=
\begin{bmatrix} 
    D_1 & A_{12} & A_{13} & ... & A_{1,n-1} & A_{1n}\\
    0 & D_2 & A_{23} & ... & A_{2,n-1} & A_{2n}\\
    0 &  0 & D_3 & ... & A_{3,n-1} & A_{3n}\\
    \vdots & \vdots & \vdots & \ddots & \vdots & \vdots\\
    0 & 0 & 0 & ... & D_{n-1} & A_{n-1,n}\\
    0 & 0 & 0 & ... & 0 & D_n      
\end{bmatrix}\in\mathcal{B}(\mathcal{H}_1\oplus \mathcal{H}_2\oplus\cdots\oplus \mathcal{H}_n),
\end{equation}
where $A:=(A_{12},\ A_{13},...,\ A_{ij},...,\ A_{n-1,n})$ consists of unknown operators $A_{ij}\in\mathcal{B}(\mathcal{H}_j,\mathcal{H}_i)$, $1\leq i<j\leq n,\ n\geq2$. Let $\mathcal{B}_n$ stand for the collection of all described tuples $A=(A_{ij})$. One easily verifies that if upper triangular operator matrix $T_n^d(A)$ is of form (\ref{OSNOVNI}), then its adjoint operator matrix $T_n^d(A)^*$ is the lower triangular operator matrix acting on $\mathcal{H}_1\oplus \mathcal{H}_2\oplus\cdots\oplus \mathcal{H}_n$. We will use this fact in duality arguments. 

Spectral properties of $T_n^d(A)$ have been widely investigated in the case $n=2$ \cite{CAO}-\cite{CAO2}, \cite{OPERATORTHEORY}-\cite{SUN}, \cite{HUANG2}-\cite{ZHANGDRUGI}, while the case of general $n$ has been neglected until a few years ago \cite{SARAJLIJA3}-\cite{WU3}.  This article is a certain continuation of the work done in the latter references.

Validity of the next lemma is clear.

\begin{Lema}\label{POMOCNALEMA}
Let $T_n^d(A)\in\mathcal{B}(X_1\oplus\cdots\oplus X_n).$ Then:
\begin{itemize}
\item[(i)] $\sigma_{SF+}(D_1)\subseteq\sigma_{SF+}(T_n^d(A))$;
\item[(ii)] $\sigma_{SF-}(D_n)\subseteq\sigma_{SF-}(T_n^d(A))$.
\end{itemize}
\end{Lema}

\begin{Teorema}(Kato's lemma, (\cite[Corollary 3.2.5]{CARADUS}))\label{KATOZATVOREN}
For $T\in\mathcal{B}(\mathcal{H})$ the following implication holds:
$$
\beta(T)<\infty\Rightarrow\mathcal{R}(T)\ is\ closed.
$$
\end{Teorema}

Finally, we list two lemmas that are crucial when it comes to using duality arguments.

\begin{Lema}\label{VEZA}
For $T\in\mathcal{B}(\mathcal{H})$ with closed range the following holds:\\
$(a)$ $\alpha(T)=\beta(T^*),\beta(T)=\alpha(T^*)$;\\
$(b)$ $T\in\Phi_+(\mathcal{H})$ if and only if $T^*\in\Phi_-(\mathcal{H})$;\\
$(c)$ $T\in\Phi_-(\mathcal{H})$ if and only if $T^*\in\Phi_+(\mathcal{H})$;\\
$(d)$ $\mathrm{ind}(T^*)=-\mathrm{ind}(T),$ if $\mathrm{ind}(T)$ is defined.
\end{Lema}

\begin{Lema}\label{DODATNALEMA}
Let $T\in\mathcal{B}(\mathcal{H})$. Then $\mathcal{R}(T)$ is closed if and only if $\mathcal{R}(T^*)$ is closed.
\end{Lema}

\section{Weylness of $T_n^d(A)$}
We begin with a result concerning upper semi-Weyl invertibility of $T_n^d(A)$.
 
\begin{Teorema}\label{LEVIVEJL}
Let $D_1\in\mathcal{B}(\mathcal{H}_1),\ D_2\in\mathcal{B}(\mathcal{H}_2),...,D_n\in\mathcal{B}(\mathcal{H}_n)$ be given. Consider the following conditions:\\[1mm]
$(i)$ $(a)$ $D_1\in\Phi_+(\mathcal{H}_1)$;\\
\hspace*{5.5mm}$(b)$ \Big($D_s\in\Phi_+(\mathcal{H}_s)$ for $2\leq s\leq n$ and $\sum\limits_{s=1}^n\alpha(D_s)\leq\sum\limits_{s=1}^n\beta(D_s)$\Big)\\ or\\ 
\hspace*{11mm}\Big($\beta(D_j)=\infty$ for some $j\in\lbrace1,...,n-1\rbrace$, $\alpha(D_s)<\infty$ for $2\leq s\leq j$ and $\mathcal{R}(D_s)$ is closed for $2\leq s\leq n$\Big);\\[3mm]
$(ii)$ There exists $A\in\mathcal{B}_n$ such that $T_n^d(A)\in\Phi_+^-(\mathcal{H}_1\oplus\cdots\oplus\mathcal{H}_n)$;\\[3mm]
$(iii)$ $(a)$ $D_1\in\Phi_+(\mathcal{H}_1)$;\\
\hspace*{8mm}$(b)$ \Big($D_s\in\Phi_+(\mathcal{H}_s)$ for $2\leq s\leq n$ and $\sum\limits_{s=1}^n\alpha(D_s)\leq\sum\limits_{s=1}^n\beta(D_s)$\Big)\\ or\\ 
\hspace*{14mm}\Big($\beta(D_j)=\infty$ for some $j\in\lbrace1,...,n-1\rbrace$, $\alpha(D_s)<\infty$ for $2\leq s\leq j$\Big).\\[3mm]
Then $(i) \Rightarrow (ii) \Rightarrow (iii)$.
\end{Teorema}
\begin{Primedba}
If $j=1$ in $(i)(b)$ or $(iii)(b)$, part "$\alpha(D_s)<\infty$ for $2\leq s\leq  j$" is omitted there.
\end{Primedba}
\begin{Primedba}
Notice the similarity between sufficent condition (i) and necessary condition (iii): parts (i)(a) and (iii)(a) are the same, while (i)(b) and (iii)(b) differ in ''$\mathcal{R}(D_s)$ is closed for $2\leq s\leq n$'' solely.
\end{Primedba}
\textbf{Proof. }$(ii)\Rightarrow (iii)$\\

\hspace*{6mm}Suppose that $T_n^d(A)$ is upper semi-Weyl. Then $T_n^d(A)$ is upper semi-Fredholm, implying $D_1\in\Phi_+(\mathcal{H}_1)$ (Lemma \ref{POMOCNALEMA}). Suppose that $(iii)(b)$ is not true. We have two possibilities. First, suppose that for $2\leq s\leq n$ we have $\beta(D_s)<\infty$. It means (Theorem \ref{KATOZATVOREN}) that $\mathcal{R}(D_s)$ is closed for $1\leq s\leq n$. Again, we have two possibilities: either there exists some $i\in\lbrace 2,...,n\rbrace$ with $\alpha(D_i)=\infty$, or we have $\sum\limits_{s=1}^n\alpha(D_s)>\sum\limits_{s=1}^n\beta(D_s)$.

Assume $\alpha(D_i)=\infty$ for some $i\in\lbrace 2,...,n\rbrace.$ We use a method from \cite{WU3}. We know that for each $A\in\mathcal{B}_n$, operator $T_n^d(A)$ regarded as an operator from $\mathcal{N}(D_1)^\bot\oplus\mathcal{N}(D_1)\oplus\mathcal{N}(D_2)^\bot\oplus\mathcal{N}(D_2)\oplus\mathcal{N}(D_3)^\bot\oplus\mathcal{N}(D_3)\oplus\cdots\oplus\mathcal{N}(D_n)^\bot\oplus\mathcal{N}(D_n)$ into $\mathcal{R}(D_1)\oplus\mathcal{R}(D_1)^\bot\oplus\mathcal{R}(D_2)\oplus\mathcal{R}(D_2)^\bot\oplus\cdots\oplus\mathcal{R}(D_{n-1})\oplus\mathcal{R}(D_{n-1})^\bot\oplus\mathcal{R}(D_n)\oplus\mathcal{R}(D_n)^\bot$ has the following block representation
\begin{equation}\label{MATRICA}
T_n^d(A)=\scalemath{0.85}
{\begin{bmatrix} 
    D_1^{(1)} & 0 & A_{12}^{(1)} & A_{12}^{(2)} & A_{13}^{(1)} & A_{13}^{(2)} & ... & A_{1n}^{(1)} & A_{1n}^{(2)}\\
    0 & 0 & A_{12}^{(3)} & A_{12}^{(4)} & A_{13}^{(3)} & A_{13}^{(4)} & ... & A_{1n}^{(3)} & A_{1n}^{(4)}\\
    0 & 0 & D_2^{(1)} & 0 & A_{23}^{(1)} & A_{23}^{(2)} & ... & A_{2n}^{(1)} & A_{2n}^{(2)}\\
    0 & 0 & 0 & 0 & A_{23}^{(3)} & A_{23}^{(4)} & ... & A_{2n}^{(3)} & A_{2n}^{(4)}\\
    0 & 0 & 0 & 0 & D_{3}^{(1)} & 0 & ... & A_{3n}^{(1)} & A_{3n}^{(2)}\\
    0 & 0 & 0 & 0 & 0 & 0 & ... & A_{3n}^{(3)} & A_{3n}^{(4)}\\
    \vdots & \vdots & \vdots & \vdots & \vdots & \vdots & \ddots & \vdots & \vdots\\
    0 & 0 & 0 & 0 & 0 & 0 & ... & A_{n-1,n}^{(1)} & A_{n-1,n}^{(2)}\\
    0 & 0 & 0 & 0 & 0 & 0 & ... & A_{n-1,n}^{(3)} & A_{n-1,n}^{(4)}\\
    0 & 0 & 0 & 0 & 0 & 0 & ... & D_n^{(1)} & 0\\
    0 & 0 & 0 & 0 & 0 & 0 & ... & 0 & 0\\
\end{bmatrix}}
\end{equation}
Evidently, $D_1^{(1)},\ D_2^{(1)},...,D_n^{(1)}$ from (\ref{MATRICA}) are invertible. Hence, there exist invertible operator matrices $U$ and $V$ so that 
\begin{equation}\label{MATRICA2}
UT_n^d(A)V=\scalemath{0.85}{\begin{bmatrix}
D_1^{(1)} & 0  & 0 & 0 & 0 & 0 & ... & 0 & 0\\
0 & 0 & 0 & A_{12}^{(4)} & 0 & A_{13}^{(4)} & ... & 0 & A_{1n}^{(4)}\\
0 & 0 & D_2^{(1)} & 0 & 0 & 0 & ... & 0 & 0\\
0 & 0 & 0 & 0 & 0 & A_{23}^{(4)} & ... & 0 & A_{2n}^{(4)}\\
0 & 0 & 0 & 0 & D_3^{(1)} & 0 & ... & 0 & 0\\
0 & 0 & 0 & 0 & 0 & 0 & ... & 0 & A_{3n}^{(4)}\\
\vdots & \vdots & \vdots & \vdots & \vdots & \vdots & \ddots & \vdots & \vdots\\
0 & 0 & 0 & 0 & 0 & 0 & ... & 0 & 0\\
0 & 0 & 0 & 0 & 0 & 0 & ... & 0 & A_{n-1,n}^{(4)}\\
0 & 0 & 0 & 0 & 0 & 0 & ... & D_n^{(1)} & 0\\
0 & 0 & 0 & 0 & 0 & 0 & ... & 0 & 0\\
\end{bmatrix}
}
\end{equation}
Operators $A_{ij}^{(4)}$ in (\ref{MATRICA}) and (\ref{MATRICA2}) are not the same, but we will keep the same notation for simplicity.
Next, it is clear that (\ref{MATRICA2}) is upper semi-Weyl if and only if 
\begin{equation}\label{MATRICA3}
\begin{bmatrix}
0 & A_{12}^{(4)} & A_{13}^{(4)} & A_{14}^{(4)} & ... & A_{1n}^{(4)}\\
0 & 0          & A_{23}^{(4)} & A_{24}^{(4)} & ... & A_{2n}^{(4)}\\
0 & 0           & 0         & A_{34}^{(4)} & ... & A_{3n}^{(4)}\\
\vdots & \vdots   &   \vdots & \vdots & \ddots & \vdots\\
0 & 0          &   0        &    0      & ... & A_{n-1,n}^{(4)}\\
0 & 0 & 0 & 0 & ... &0
\end{bmatrix}
:
\begin{bmatrix}
\mathcal{N}(D_1)\\
\mathcal{N}(D_2)\\
\mathcal{N}(D_3)\\
\mathcal{N}(D_4)\\
\vdots\\
\mathcal{N}(D_n)
\end{bmatrix}
\rightarrow
\begin{bmatrix}
\mathcal{R}(D_1)^\bot\\
\mathcal{R}(D_2)^\bot\\
\mathcal{R}(D_3)^\bot\\
\vdots\\
\mathcal{R}(D_{n-1})^\bot\\
\mathcal{R}(D_{n})^\bot
\end{bmatrix}
\end{equation}
is upper semi-Weyl. Since $\sum\limits_{s=1}^{i-1}\beta(D_s)<\infty$ and $\alpha(D_i)=\infty$, it follows that 
$$
\alpha\left(\begin{bmatrix}
A_{1i}^{(4)}\\
A_{2i}^{(4)}\\
A_{3i}^{(4)}\\
\vdots\\
A_{i-1,i}^{(4)}
\end{bmatrix}\right)=\infty,
$$
and hence operator defined in (\ref{MATRICA3}) is not upper semi-Weyl for every $A\in\mathcal{B}_n$. This proves the desired.

Assume next that $\alpha(D_s)<\infty$ for $2\leq s\leq n$. Then $\sum\limits_{s=1}^n\alpha(D_s)>\sum\limits_{s=1}^n\beta(D_s)$, and for each $A\in\mathcal{B}_n$, $T_n^d(A)$ has representation as (\ref{MATRICA}), and we use (\ref{MATRICA2}) and (\ref{MATRICA3}) again. Since $D_s$, $1\leq s\leq n$ are upper semi-Fredholm, then $T_n^d(A)$ is upper semi-Weyl if and only if (\ref{MATRICA3}) is upper semi-Weyl. But $\sum\limits_{s=1}^n\beta(D_s)<\sum\limits_{s=1}^n\alpha(D_s)$ implies (\ref{MATRICA3}) is not upper semi-Weyl for every $A\in\mathcal{B}_n$. Contradiction.

Second option is that there is $j\in\lbrace 2
,...,n\rbrace$ with $\beta(D_j)=\infty$, and assume we have found the smallest such $j$. Then $\beta(D_s)<\infty$ for $1\leq s\leq j-1$, hence $\mathcal{R}(D_s)$ is closed for $1\leq s\leq j-1$. Now, $\alpha(D_s)<\infty$ for $2\leq s\leq j-1$ is not possible, otherwise $(iii)(b)$ would be true. Finally, $\alpha(D_j)=\infty$ for some $j\in\lbrace2,...,j-1\rbrace$ and be proceed with (\ref{MATRICA}), (\ref{MATRICA2}), (\ref{MATRICA3}) applied to $T_{j-1}^d(A)$.

\noindent$(i)\Rightarrow(ii)$

Assume that $D_1\in\Phi_{+}(\mathcal{H}_1)$ and $(i)(b)$ holds. If $D_s\in\Phi_{+
}(\mathcal{H}_j)$ for $2\leq s\leq n$ and $\sum\limits_{s=1}^n\alpha(D_s)\leq\sum\limits_{s=1}^n\beta(D_s)$, we choose trivially $A=\mathbf{0}$ and $T_n^d(A)$ is upper semi-Weyl. \\
Suppose that $\beta(D_j)=\infty$ for some $j\in\lbrace 1,...,n-1\rbrace,$ $\alpha(D_s)<\infty$ for $2\leq s\leq j$ and $\mathcal{R}(D_s)$ is closed for all $1\leq s\leq n$. Assume that $\lbrace f_s^{(k)}\rbrace_{s=1}^\infty$, $\lbrace e_s^{(1)}\rbrace_{s=1}^\infty$,  $\lbrace e_s^{(2)}\rbrace_{s=1}^\infty$,..., $\lbrace e_s^{(n-1)}\rbrace_{s=1}^\infty$ are orthogonal bases of $\mathcal{R}(D_k)^\bot$, $\mathcal{H}_2$,..., $\mathcal{H}_n$, respectively. We have two cases. Again, we adopt a method from \cite{WU3}.

\textbf{Case 1}: $\beta(D_1)=\infty$

In this case it holds  $\alpha(D_1)<\infty$, $\mathcal{R}(D_s)$ is closed for all $1\leq s\leq n$ and $\beta(D_1)=\infty$. We find $A\in\mathcal{B}_n$ such that $\alpha(T_n^d(A))<\infty$ and $\mathcal{R}(T_n^d(A))$ is closed.  We choose $A=(A_{ij})_{1\leq i<j\leq n}$ so that $A_{ij}=0$ if $i>1$, that is we place all nonzero operators of tuple $A$ in the first row. It remains to define $A_{1s}$ for $2\leq s\leq n$. We put
$$
\begin{aligned}
A_{12}(e_s^{(1)})=f_{ns}^{(1)},\quad s=1,2,...;\\
A_{13}(e_s^{(2)})=f_{ns+1}^{(1)},\quad s=1,2,...;\\
\cdots\\
A_{1n}(e_s^{(n-1)})=f_{ns+n-2}^{(1)},\quad s=1,2,....
\end{aligned}
$$

Now we have chosen our $A=(A_{ij})$, it is easy to show that $\mathcal{N}(T_n^d(A))=\mathcal{N}(D_1)\oplus\lbrace\mathbf{0}\rbrace\oplus\cdots\oplus\lbrace\mathbf{0}\rbrace$. Therefore, $\alpha(T_n^d(A))=\alpha(D_1)<\infty$.

Secondly, we show that $\mathcal{R}(T_n^d(A))$ is closed and $\beta(T_n^d(A))=\infty$. Since $\mathcal{R}(D_s)$ is closed for all $1\leq s\leq n$, it will follow that $\mathcal{R}(T_n^d(A))$ is closed if we prove that $\mathcal{R}(A_{1s})$ is closed for $2\leq s\leq n$. But, since we are in the setting of separable Hilbert spaces, with regards to definition of $A_{1s}$'s, the former is obvious. We have that $T_n^d(A)$ is upper semi-Fredholm, and since $\beta(T_n^d(A))=\beta(D_1)=\infty$ due to definition of $A_{1s}$'s, we find that $T_n^d(A)$ is upper semi-Weyl.\\

\noindent\textbf{Case 2}: $\beta(D_k)=\infty$ for some $k\in\lbrace2,...,n-1\rbrace$

In this case it holds  $\alpha(D_s)<\infty$, $1\leq s\leq k$, $\mathcal{R}(D_s)$ is closed for all $1\leq s\leq n$ and $\beta(D_k)=\infty$. We find $A\in\mathcal{B}_n$ such that $\alpha(T_n^d(A))<\infty$ and $\mathcal{R}(T_n^d(A))$ is closed.  We choose $A=(A_{ij})_{1\leq i<j\leq n}$ so that $A_{ij}=0$ if $i\neq k$, that is we place all nonzero operators of tuple $A$ in the $k$-th row. It remains to define $A_{ks}$ for $k+1\leq s\leq n$. We put
$$
\begin{aligned}
A_{k,k+1}(e_s^{(k)})=f_{ns}^{(k)},\quad s=1,2,...;\\
A_{k,k+2}(e_s^{(k+1)})=f_{ns+1}^{(k)},\quad s=1,2,...;\\
\cdots\\
A_{kn}(e_s^{(n-1)})=f_{ns+n-k-1}^{(k)},\quad s=1,2,....
\end{aligned}
$$

Now we have chosen our $A=(A_{ij})$
, it is easy to show that $\mathcal{N}(T_n^d(A))=\mathcal{N}(D_1)\oplus\cdots\oplus\mathcal{N}(D_{k})\oplus\lbrace\mathbf{0}\rbrace\oplus\cdots\oplus\lbrace\mathbf{0}\rbrace$. Therefore, $\alpha(T_n^d(A))\leq\alpha(D_1)+\cdots+\alpha(D_{k})<\infty$.

Secondly, we show that $\mathcal{R}(T_n^d(A))$ is closed and $\beta(T_n^d(A))=\infty$. Since $\mathcal{R}(D_s)$ is closed for all $1\leq s\leq n$, it will follow that $\mathcal{R}(T_n^d(A))$ is closed if we prove that $\mathcal{R}(A_{ks})$ is closed for $k+1\leq s\leq n$. But, since we are in the setting of separable Hilbert spaces, with regards to definition of $A_{ks}$'s, the former is obvious. We have that $T_n^d(A)$ is upper semi-Fredholm, and since $\beta(T_n^d(A))=\beta(D_k)=\infty$ due to definition of $A_{ks}$'s, we find that $T_n^d(A)$ is upper semi-Weyl.\\
 $\square$
\begin{Primedba}\label{VAZISVUDA}
Notice the validity of part $(ii)\Rightarrow(iii)$ without assuming separability of $\mathcal{H}_1,...,\mathcal{H}_n$.
\end{Primedba}

Next corollary is immediate from Theorem \ref{LEVIVEJL}. 
\begin{Posledica}(\cite[Theorem 2.5]{WU3}, corrected version)\label{POSLEDICA}
Let $D_1\in\mathcal{B}(\mathcal{H}_1),\ D_2\in\mathcal{B}(\mathcal{H}_2),...,D_n\in\mathcal{B}(\mathcal{H}_n)$. Then
$$
\begin{aligned}
\sigma_{SF+}(D_1)\cup\Big(\bigcup\limits_{k=2}^{n+1}\Delta_k\Big)\subseteq\\\bigcap\limits_{A\in\mathcal{B}_n}\sigma_{aw}(T_n^d(A))\subseteq\\\sigma_{SF+}(D_1)\cup\Big(\bigcup\limits_{k=2}^{n+1}\Delta_k\Big)\cup\Big(\bigcup\limits_{k=2}^n\Delta_k'\Big),
\end{aligned}
$$
where
$$
\Delta_k:=\Big\lbrace\lambda\in\mathds{C}:\ \alpha(D_k-\lambda)=\infty\ and\  \sum\limits_{s=1}^{k-1}\beta(D_s-\lambda)<\infty\Big\rbrace,\ 2\leq k\leq n,
$$
$$
\Delta_{n+1}:=\Big\lbrace\lambda\in\mathds{C}:\ \sum\limits_{s=1}^n\beta(D_s-\lambda)<\sum\limits_{s=1}^n\alpha(D_s-\lambda)\Big\rbrace,
$$
$$
\Delta_k':=\Big\lbrace\lambda\in\mathds{C}:\ \mathcal{R}(D_k-\lambda)\ is\ not\ closed\Big\rbrace,\ 2\leq k\leq n.
$$
\end{Posledica}
\begin{Primedba}
One should also spot a difference between collections $\Delta_k$, $2\leq k\leq n$, in Corollary \ref{POSLEDICA} and in \cite[Theorem 2.5]{WU3}. This difference is implied by the existence of sets $\Delta_k'$, $2\leq k\leq n$, in the formulation of Corollary \ref{POSLEDICA}. Our estimates are better in a sense that $\Delta_k$ in Corollary \ref{POSLEDICA} is a subset of $\Delta_k$ from \cite{WU3} for every $2\leq k\leq n$.
\end{Primedba}

Previous statements for $n=2$ become very simple, as shown in the sequel.
\begin{Teorema}(\cite[Theorem 2.5]{ZHANG}, corrected version)\label{LEVIVEJLn=2}
Let $D_1\in\mathcal{B}(\mathcal{H}_1), D_2\in\mathcal{B}(\mathcal{H}_2)$. Consider the following statements:\\[1mm]
$(i)$ $(a)$ $D_1\in\Phi_+(\mathcal{H}_1)$;\\
\hspace*{5.5mm}$(b)$  \Big($D_2\in\Phi_+(\mathcal{H}_2)$ and $\alpha(D_1)+\alpha(D_2)\leq\beta(D_1)+\beta(D_2)$\Big) \\or\\ \hspace*{11mm}\Big($\beta(D_1)=\infty$ and $\mathcal{R}(D_2)$ is closed\Big);\\[1mm]
$(ii)$ There exists $A\in\mathcal{B}_2$ such that $T_2^d(A)\in\Phi_+^-(\mathcal{H}_1\oplus\mathcal{H}_2)$;\\[1mm]
$(iii)$ $(a)$ $D_1\in\Phi_+(\mathcal{H}_1)$;\\
\hspace*{8.2mm}$(b)$ \Big($D_2\in\Phi_+(\mathcal{H}_2)$ and $\alpha(D_1)+\alpha(D_2)\leq\beta(D_1)+\beta(D_2)$\Big) \\or\\ \hspace*{13mm}$\beta(D_1)=\infty$.\\[1mm]
Then $(i)\Rightarrow(ii)\Rightarrow(iii)$.
\end{Teorema}

Notice that Theorem \ref{LEVIVEJLn=2} is a corrected version of \cite[Theorem 2.5]{ZHANG}. Condition '$\mathcal{R}(D_2)$ is closed' in $(i)(b)$ is omitted in \cite{ZHANG}, which is an oversight. Without that condition we can not prove that $\mathcal{R}(T_2^d(A))$ is closed and therefore direction $(ii)\Rightarrow(i)$ in \cite[Theorem 2.5]{ZHANG} would not hold. 
\begin{Posledica}(\cite[Corollary 2.7]{ZHANG}, corrected version)\label{POSLEDICA2}
Let $D_1\in\mathcal{B}(\mathcal{H}_1),\ D_2\in\mathcal{B}(\mathcal{H}_2)$. Then
$$
\sigma_{SF+}(D_1)\cup\Delta\cup\Delta'\subseteq\bigcap\limits_{A\in\mathcal{B}_2}\sigma_{aw}(T_2^d(A))\subseteq\sigma_{SF+}(D_1)\cup\Delta\cup\Delta'\cup\Delta'',
$$
where
$$
\Delta:=\Big\lbrace\lambda\in\mathds{C}:\ \alpha(D_2-\lambda)=\infty\ and\  \beta(D_1-\lambda)<\infty\Big\rbrace,
$$
$$
\Delta':=\Big\lbrace\lambda\in\mathds{C}:\ \beta(D_1-\lambda)+\beta(D_2-\lambda)<\alpha(D_1-\lambda)+\alpha(D_2-\lambda)\Big\rbrace,
$$
$$
\Delta'':=\lbrace\lambda\in\mathds{C}:\ \mathcal{R}(D_2-\lambda)\ is\ not\ closed\rbrace.
$$

\end{Posledica}

Now we list statements dealing with the lower semi-Weyl spectrum.
\begin{Teorema}\label{DESNIVEJL}
Let $D_1\in\mathcal{B}(\mathcal{H}_1),\ D_2\in\mathcal{B}(\mathcal{H}_2),...,D_n\in\mathcal{B}(\mathcal{H}_n)$ be given. Consider the following conditions:\\
$(i)$   $(a)$ $D_n\in\Phi_-(\mathcal{H}_n)$;\\
\hspace*{5.5mm}$(b)$ \Big($D_s\in\Phi_-(\mathcal{H}_s)$ for $1\leq s\leq n-1$ and $\sum\limits_{s=1}^n\beta(D_s)\leq\sum\limits_{s=1}^n\alpha(D_s)$\Big) \\or\\ \hspace*{11mm}\Big($\alpha(D_j)=\infty$ for some $j\in\lbrace2,...,n\rbrace$, $\beta(D_s)<\infty$ for $j\leq s\leq n-1$ and $\mathcal{R}(D_s)$ is closed for $1\leq s\leq n-1$\Big);\\[3mm]
$(ii)$ There exists $A\in\mathcal{B}_n$ such that $T_n^d(A)\in\Phi_-^+(\mathcal{H}_1\oplus\cdots\oplus\mathcal{H}_n)$;\\[3mm]
$(iii)$ $(a)$ $D_n\in\Phi_-(\mathcal{H}_n)$;\\
\hspace*{8.2mm}$(b)$ \Big($D_s\in\Phi_-(\mathcal{H}_s)$ for $1\leq s\leq n-1$ and $\sum\limits_{s=1}^n\beta(D_s)\leq\sum\limits_{s=1}^n\alpha(D_s)$\Big)\\ or \\\hspace*{14mm}\Big($\alpha(D_j)=\infty$ for some $j\in\lbrace2,...,n\rbrace$, $\beta(D_s)<\infty$ for $j\leq s\leq n-1$\Big).\\[3mm]
Then $(i) \Rightarrow (ii) \Rightarrow (iii)$.
\end{Teorema}
\begin{Primedba}
If $j=n$ in $(i)(b)$ or $(iii)(b)$, part ''$\beta(D_s)<\infty$ for $j\leq s\leq n-1$'' is omitted there.
\end{Primedba}
\begin{Primedba}
Notice the similarity between sufficent condition (i) and necessary condition (iii): parts (i)(a) and (iii)(a) are the same, while (i)(b) and (iii)(b) differ in ''$\mathcal{R}(D_s)$ is closed for $1\leq s\leq n-1$'' solely.
\end{Primedba}
\textbf{Proof. }This easily follows from the statement of Theorem \ref{LEVIVEJL} by duality argument, putting into use Lemmas \ref{VEZA} and \ref{DODATNALEMA}. $\square$

\begin{Posledica}(\cite[Theorem 2.6]{WU3}, corrected version)\label{POSLEDICA3}
Let $D_1\in\mathcal{B}(\mathcal{H}_1),\ D_2\in\mathcal{B}(\mathcal{H}_2),...,D_n\in\mathcal{B}(\mathcal{H}_n)$. Then
$$
\begin{aligned}
\sigma_{SF-}(D_n)\cup\Big(\bigcup\limits_{k=1}^{n-1}\Delta_k\Big)\cup\Delta_{n+1}\subseteq\\\bigcap\limits_{A\in\mathcal{B}_n}\sigma_{sw}(T_n^d(A))\subseteq\\\sigma_{SF-}(D_n)\cup\Big(\bigcup\limits_{k=1}^{n-1}\Delta_k\Big)\cup\Delta_{n+1}\cup\Big(\bigcup\limits_{k=1}^{n-1}\Delta_k'\Big),
\end{aligned}
$$
where
$$
\Delta_k:=\Big\lbrace\lambda\in\mathds{C}:\ \beta(D_k-\lambda)=\infty\ and\   \sum\limits_{s=k+1}^{n}\alpha(D_s-\lambda)<\infty\Big\rbrace,\ 1\leq k\leq n-1,
$$
$$
\Delta_{n+1}:=\Big\lbrace\lambda\in\mathds{C}:\ \sum\limits_{s=1}^n\alpha(D_s-\lambda)<\sum\limits_{s=1}^n\beta(D_s-\lambda)\Big\rbrace,
$$
$$
\Delta_k':=\Big\lbrace\lambda\in\mathds{C}:\ \mathcal{R}(D_k-\lambda)\ is\ not\ closed\Big\rbrace,\ 1\leq k\leq n-1.
$$

\end{Posledica}

\begin{Primedba}
Again, existence of sets $\Delta_k'$, $1\leq k\leq n-1$ in the statement of Corollary \ref{POSLEDICA3} implies a difference between definitions of collections $\Delta_k$, $1\leq k\leq n-1$ in Corollary \ref{POSLEDICA3} and in \cite[Theorem 2.6]{WU3}.  
\end{Primedba}

If we put $n=2$ we get:

\begin{Teorema}(\cite[Theorem 2.6]{ZHANG}, corrected version)\label{DESNIVEJLn=2}
Let $D_1\in\mathcal{B}(\mathcal{H}_1), D_2\in\mathcal{B}(\mathcal{H}_2)$. Consider the following conditions:\\[1mm]
$(i)$ $(a)$ $D_2\in\Phi_-(\mathcal{H}_2)$;\\
\hspace*{5.5mm}$(b)$ \Big($D_1\in\Phi_-(\mathcal{H}_1)$ and $\alpha(D_1)+\alpha(D_2)\geq\beta(D_1)+\beta(D_2)$\Big) \\or\\ \hspace*{11mm}\Big($\alpha(D_2)=\infty$ and $\mathcal{R}(D_1)$ is closed\Big);\\[1mm]
$(ii)$ There exists $A\in\mathcal{B}_2$ such that $T_2^d(A)\in\Phi_-^+(\mathcal{H}_1\oplus\mathcal{H}_2)$;\\[1mm]
$(iii)$ $(a)$ $D_2\in\Phi_-(\mathcal{H}_2)$;\\
\hspace*{8.2mm}$(c)$ \Big($D_1\in\Phi_-(\mathcal{H}_1)$ and $\alpha(D_1)+\alpha(D_2)\geq\beta(D_1)+\beta(D_2)$\Big) \\or\\ \hspace*{13mm} $\alpha(D_2)=\infty$.\\[3mm]
Then $(i)\Rightarrow(ii)\Rightarrow(iii)$.
\end{Teorema}
\begin{Posledica}(\cite[Corollary 2.8]{ZHANG}, corrected version)\label{POSLEDICA4}
Let $D_1\in\mathcal{B}(\mathcal{H}_1),\ D_2\in\mathcal{B}(\mathcal{H}_2)$. Then
$$
\sigma_{SF-}(D_2)\cup\Delta\cup\Delta'\subseteq\bigcap\limits_{A\in\mathcal{B}_2}\sigma_{sw}(T_2^d(A))\subseteq\sigma_{SF-}(D_2)\cup\Delta\cup\Delta'\cup\Delta'',
$$
where
$$
\Delta:=\Big\lbrace\lambda\in\mathds{C}:\ \beta(D_1-\lambda)=\infty\ and\  \alpha(D_2-\lambda)<\infty\Big\rbrace,
$$
$$
\Delta':=\Big\lbrace\lambda\in\mathds{C}:\ \alpha(D_1-\lambda)+\alpha(D_2-\lambda)<\beta(D_1-\lambda)+\beta(D_2-\lambda)\Big\rbrace,
$$
$$
\Delta'':=\lbrace\lambda\in\mathds{C}:\ \mathcal{R}(D_1-\lambda)\ is\ not\  closed\rbrace.
$$

\end{Posledica}
\subsection{Fredholmness of $T_n^d(A)$}
In this section we provide statements related to the Fredholmness of $T_n^d(A)$. One can notice that these statements are quite similar to the ones presented in the previous section. Their proofs are also similar, and so we omit them.

We start with a result which deals with the upper semi-Fredholm invertibility of $T_n^d(A)$. 
\begin{Teorema}\label{LEVIFREDHOLM}
Let $D_1\in\mathcal{B}(\mathcal{H}_1),\ D_2\in\mathcal{B}(\mathcal{H}_2),...,D_n\in\mathcal{B}(\mathcal{H}_n)$ be given. Consider the following conditions:\\[1mm]
$(i)$ $(a)$ $D_1\in\Phi_+(\mathcal{H}_1)$;\\
\hspace*{5.5mm}$(b)$ $D_s\in\Phi_+(\mathcal{H}_s)$ for $2\leq s\leq n$\\ or\\ 
\hspace*{11mm}\Big($\beta(D_j)=\infty$ for some $j\in\lbrace1,...,n-1\rbrace$, $\alpha(D_s)<\infty$ for $2\leq s\leq j$ and $\mathcal{R}(D_s)$ is closed for $2\leq s\leq n$\Big);\\[3mm]
$(ii)$ There exists $A\in\mathcal{B}_n$ such that $T_n^d(A)\in\Phi_+(\mathcal{H}_1\oplus\cdots\oplus\mathcal{H}_n)$;\\[3mm]
$(iii)$ $(a)$ $D_1\in\Phi_+(\mathcal{H}_1)$;\\
\hspace*{8mm}$(b)$ $D_s\in\Phi_+(\mathcal{H}_s)$ for $2\leq s\leq n$\\ or\\ 
\hspace*{14mm}\Big($\beta(D_j)=\infty$ for some $j\in\lbrace1,...,n-1\rbrace$, $\alpha(D_s)<\infty$ for $2\leq s\leq j$\Big).\\[3mm]
Then $(i) \Rightarrow (ii) \Rightarrow (iii)$.
\end{Teorema}
\begin{Primedba}
If $j=1$ in $(i)(b)$ or $(iii)(b)$, part "$\alpha(D_s)<\infty$ for $2\leq s\leq  j$" is omitted there.
\end{Primedba}
\begin{Primedba}
Notice the similarity between sufficent condition (i) and necessary condition (iii): parts (i)(a) and (iii)(a) are the same, while (i)(b) and (iii)(b) differ in ''$\mathcal{R}(D_s)$ is closed for $2\leq s\leq n$'' solely.
\end{Primedba}
\begin{Primedba}\label{VAZISVUDA2}
Again, we have the validity of part $(ii)\Rightarrow(iii)$ without assuming separability of $\mathcal{H}_1,...,\mathcal{H}_n$.
\end{Primedba}

\begin{Posledica}(\cite[Theorem 1]{WU2}, corrected version)\label{POSLEDICA5}
Let $D_1\in\mathcal{B}(\mathcal{H}_1),\ D_2\in\mathcal{B}(\mathcal{H}_2),...,D_n\in\mathcal{B}(\mathcal{H}_n)$. Then
$$
\begin{aligned}
\sigma_{SF+}(D_1)\cup\Big(\bigcup\limits_{k=2}^{n}\Delta_k\Big)\subseteq\\\bigcap\limits_{A\in\mathcal{B}_n}\sigma_{SF+}(T_n^d(A))\subseteq\\\sigma_{SF+}(D_1)\cup\Big(\bigcup\limits_{k=2}^{n}\Delta_k\Big)\cup\Big(\bigcup\limits_{k=2}^n\Delta_k'\Big),
\end{aligned}
$$
where
$$
\Delta_k:=\Big\lbrace\lambda\in\mathds{C}:\ \alpha(D_k-\lambda)=\infty\ and\  \sum\limits_{s=1}^{k-1}\beta(D_s-\lambda)<\infty\Big\rbrace,\ 2\leq k\leq n,
$$
$$
\Delta_k':=\Big\lbrace\lambda\in\mathds{C}:\ \mathcal{R}(D_k-\lambda)\ is\ not\ closed\Big\rbrace,\ 2\leq k\leq n.
$$
\end{Posledica}
\begin{Primedba}
Notice  a difference between definitions of sets $\Delta_k$, $2\leq k\leq n$, in Corollary \ref{POSLEDICA5} and in \cite[Theorem 1]{WU2}.  
\end{Primedba}
If we put $n=2$ we get:
\begin{Teorema}(\cite[Theorem 2.10]{ZHANG}, corrected version)\label{LEVIFREDHOLMn=2}
Let $D_1\in\mathcal{B}(\mathcal{H}_1), D_2\in\mathcal{B}(\mathcal{H}_2)$. Consider the following statements:\\[1mm]
$(i)$ $(a)$ $D_1\in\Phi_+(\mathcal{H}_1)$;\\
\hspace*{5.5mm}$(b)$  $D_2\in\Phi_+(\mathcal{H}_2)$ or \Big($\beta(D_1)=\infty$ and $\mathcal{R}(D_2)$ is closed\Big).\\[1mm]
$(ii)$ There exists $A\in\mathcal{B}_2$ such that $T_2^d(A)\in\Phi_+(\mathcal{H}_1\oplus\mathcal{H}_2)$;\\[1mm]
$(iii)$ $(a)$ $D_1\in\Phi_+(\mathcal{H}_1)$;\\
\hspace*{8.2mm}$(b)$ $D_2\in\Phi_+(\mathcal{H}_2)$ or $\beta(D_1)=\infty$.\\[1mm]
Then $(i)\Rightarrow(ii)\Rightarrow(iii)$.
\end{Teorema}

\begin{Posledica}(\cite[Corollary 2.12]{ZHANG}, corrected version)\label{POSLEDICA6}
Let $D_1\in\mathcal{B}(\mathcal{H}_1),\ D_2\in\mathcal{B}(\mathcal{H}_2)$. Then
$$
\sigma_{SF+}(D_1)\cup\Delta\subseteq\bigcap\limits_{A\in\mathcal{B}_2}\sigma_{SF+}(T_2^d(A))\subseteq\sigma_{SF+}(D_1)\cup\Delta\cup\Delta',
$$
where
$$
\Delta:=\Big\lbrace\lambda\in\mathds{C}:\ \alpha(D_2-\lambda)=\infty\ and\  \beta(D_1-\lambda)<\infty\Big\rbrace,
$$
$$
\Delta':=\lbrace\lambda\in\mathds{C}:\ \mathcal{R}(D_2-\lambda)\ is\ not\ closed\rbrace.
$$

\end{Posledica}

Now we list statements dealing with the lower semi-Fredholm spectrum.
\begin{Teorema}\label{DESNIFREDHOLM}
Let $D_1\in\mathcal{B}(\mathcal{H}_1),\ D_2\in\mathcal{B}(\mathcal{H}_2),...,D_n\in\mathcal{B}(\mathcal{H}_n)$ be given. Consider the following conditions:\\
$(i)$   $(a)$ $D_n\in\Phi_-(\mathcal{H}_n)$;\\
\hspace*{5.5mm}$(b)$ $D_s\in\Phi_-(\mathcal{H}_s)$ for $1\leq s\leq n-1$ \\or\\ \hspace*{11mm}\Big($\alpha(D_j)=\infty$ for some $j\in\lbrace2,...,n\rbrace$, $\beta(D_s)<\infty$ for $j\leq s\leq n-1$ and $\mathcal{R}(D_s)$ is closed for $1\leq s\leq n-1$\Big);\\[3mm]
$(ii)$ There exists $A\in\mathcal{B}_n$ such that $T_n^d(A)\in\Phi_-(\mathcal{H}_1\oplus\cdots\oplus\mathcal{H}_n)$;\\[3mm]
$(iii)$ $(a)$ $D_n\in\Phi_-(\mathcal{H}_n)$;\\
\hspace*{8.2mm}$(b)$ $D_s\in\Phi_-(\mathcal{H}_s)$ for $1\leq s\leq n-1$\\ or \\\hspace*{14mm}\Big($\alpha(D_j)=\infty$ for some $j\in\lbrace2,...,n\rbrace$, $\beta(D_s)<\infty$ for $j\leq s\leq n-1$\Big).\\[3mm]
Then $(i) \Rightarrow (ii) \Rightarrow (iii)$.
\end{Teorema}
\begin{Primedba}
If $j=n$ in $(i)(b)$ or $(iii)(b)$, part ''$\beta(D_s)<\infty$ for $j\leq s\leq n-1$'' is omitted there.
\end{Primedba}
\begin{Primedba}
Notice the similarity between sufficent condition (i) and necessary condition (iii): parts (i)(a) and (iii)(a) are the same, while (i)(b) and (iii)(b) differ in ''$\mathcal{R}(D_s)$ is closed for $1\leq s\leq n-1$'' solely.
\end{Primedba}

\begin{Posledica}(\cite[Theorem 2]{WU2}, corrected version)\label{POSLEDICA7}
$$
\begin{aligned}
\sigma_{SF-}(D_n)\cup\Big(\bigcup\limits_{k=1}^{n-1}\Delta_k\Big)\subseteq\\\bigcap\limits_{A\in\mathcal{B}_n}\sigma_{SF-}(T_n^d(A))\subseteq\\\sigma_{SF-}(D_n)\cup\Big(\bigcup\limits_{k=1}^{n-1}\Delta_k\Big)\cup\Big(\bigcup\limits_{k=1}^{n-1}\Delta_k'\Big),
\end{aligned}
$$
where
$$
\Delta_k:=\Big\lbrace\lambda\in\mathds{C}:\ \beta(D_k-\lambda)=\infty\ and\   \sum\limits_{s=k+1}^{n}\alpha(D_s-\lambda)<\infty\Big\rbrace,\ 1\leq k\leq n-1,
$$
$$
\Delta_k':=\Big\lbrace\lambda\in\mathds{C}:\ \mathcal{R}(D_k-\lambda)\ is\ not\ closed\Big\rbrace,\ 1\leq k\leq n-1.
$$
\end{Posledica}
\begin{Primedba}
Again we have a difference between definitions of the sets $\Delta_k$, $1\leq k\leq n-1$ in Corollary \ref{POSLEDICA7} and in \cite[Theorem 2]{WU2}.  
\end{Primedba}
\begin{Teorema}(\cite[Theorem 2.11]{ZHANG}, corrected version)\label{DESNIFREDHOLMn=2}
Let $D_1\in\mathcal{B}(\mathcal{H}_1), D_2\in\mathcal{B}(\mathcal{H}_2)$. Consider the following conditions:\\[1mm]
$(i)$ $(a)$ $D_2\in\Phi_-(\mathcal{H}_2)$;\\
\hspace*{5.5mm}$(b)$ $D_1\in\Phi_-(\mathcal{H}_1)$ or \Big($\alpha(D_2)=\infty$ and $\mathcal{R}(D_1)$ is closed\Big);\\[1mm]
$(ii)$ There exists $A\in\mathcal{B}_2$ such that $T_2^d(A)\in\Phi_-(\mathcal{H}_1\oplus\mathcal{H}_2)$;\\[1mm]
$(iii)$ $(a)$ $D_2\in\Phi_-(\mathcal{H}_2)$;\\
\hspace*{8.2mm}$(c)$ $D_1\in\Phi_-(\mathcal{H}_1)$ or $\alpha(D_2)=\infty$.\\[3mm]
Then $(i)\Rightarrow(ii)\Rightarrow(iii)$.
\end{Teorema}
\begin{Posledica}(\cite[Corollary 2.13]{ZHANG}, corrected version)\label{POSLEDICA8}
Let $D_1\in\mathcal{B}(\mathcal{H}_1),\ D_2\in\mathcal{B}(\mathcal{H}_2)$. Then
$$
\sigma_{SF-}(D_2)\cup\Delta\subseteq\bigcap\limits_{A\in\mathcal{B}_2}\sigma_{SF-}(T_2^d(A))\subseteq\sigma_{SF-}(D_2)\cup\Delta\cup\Delta',
$$
where
$$
\Delta:=\Big\lbrace\lambda\in\mathds{C}:\ \beta(D_1-\lambda)=\infty\ and\  \alpha(D_2-\lambda)<\infty\Big\rbrace,
$$
$$
\Delta':=\lbrace\lambda\in\mathds{C}:\ \mathcal{R}(D_1-\lambda)\ is\ not\  closed\rbrace.
$$
\end{Posledica}

And this is a result about Fredholm invertibility of $T_n^d(A)$.
\begin{Teorema}\label{FREDHOLM}
Let $D_1\in\mathcal{B}(\mathcal{H}_1),\ D_2\in\mathcal{B}(\mathcal{H}_2),...,D_n\in\mathcal{B}(\mathcal{H}_n)$. Consider the following statements:\\
$(i)$   $(a)$ $D_1\in\Phi_+(\mathcal{H}_1)$ and $D_n\in\Phi_-(\mathcal{H}_n)$;\\
\hspace*{5.5mm}$(b)$ \Big($D_j\in\Phi_+(\mathcal{H}_j)$ for $2\leq j\leq n$ and $D_k\in\Phi_-(\mathcal{H}_k)$ for $1\leq k\leq n-1$\Big)\\ or\\\hspace*{11mm} \Big($\beta(D_j)=\infty$ for some $j\in\lbrace1,...,n-1\rbrace$, $\alpha(D_j)<\infty$, $\alpha(D_k)=\infty$ for some $k\in\lbrace2,...,n\rbrace$, $k>j$, $\beta(D_k)<\infty$, $\alpha(D_s),\beta(D_s)<\infty$ for $1\leq s\leq j-1$ and $k+1\leq s\leq n$, and $\mathcal{R}(D_s)$ is closed for $2\leq s\leq n-1$\Big)\\[3mm]
$(ii)$ There exists $A\in\mathcal{B}_n$ such that $T_n^d(A)\in\Phi(\mathcal{H}_1\oplus\cdots\oplus\mathcal{H}_n)$;\\[3mm]
$(iii)$ $(a)$ $D_1\in\Phi_+(\mathcal{H}_1)$ and  $D_n\in\Phi_-(\mathcal{H}_n)$;\\
\hspace*{8mm}$(b)$ \Big($D_j\in\Phi_+(\mathcal{H}_j)$ for $2\leq j\leq n$ and $D_k\in\Phi_-(\mathcal{H}_k)$ for $1\leq k\leq n-1$\Big)\\ or\\ \hspace*{13mm}\Big($\beta(D_j)=\infty$ for some $j\in\lbrace1,...,n-1\rbrace$ and $\alpha(D_s)<\infty$ for $2\leq s\leq j$, $\alpha(D_k)=\infty$ for some $k\in\lbrace2,...,n\rbrace$, and $\beta(D_s)<\infty$ for $k\leq s\leq n-1$, $k>j$\Big).\\[3mm]
Then $(i) \Rightarrow (ii) \Rightarrow (iii)$.
\end{Teorema}
\begin{Primedba}
If $j=1$ and/or $k=n$ in $(i)(b)$, condition that is ought to hold for $1\leq s\leq j-1$ and/or $k+1\leq s\leq n$ is omitted there.\\[1mm]
If $j=1$ and/or $k=n$ in $(iii)(b)$, condition that is ought to hold for $2\leq s\leq j$ and/or $k\leq s\leq n-1$ is omitted there.\\[1mm]
If $n=2$, condition "$\mathcal{R}(D_s)$ is closed for $2\leq s\leq n-1$" is omitted in $(i)(b)$.
\end{Primedba}
\begin{Primedba}
Notice the similarity between sufficent condition (i) and necessary condition (iii): again, parts (i)(a) and (iii)(a) are the same, while (i)(b) and (iii)(b) differ only slightly.
\end{Primedba}
\textbf{Proof. }$(ii)\Rightarrow(iii)$

Let $T_n^d(A)$ be Fredholm for some $A\in\mathcal{B}_n$. Then $T_n^d(A)$ is both upper and lower semi-Fredholm, and so by employing Theorems \ref{LEVIFREDHOLM} and \ref{DESNIFREDHOLM} we easily get the desired.

\noindent$(i)\Rightarrow(ii)$

Let the conditions $(i)(a)$ and $(i)(b)$ hold. If $D_j\in\Phi_+(\mathcal{H}_j)$ for $2\leq j\leq n$ and $D_k\in\Phi_-(\mathcal{H}_k)$ for $1\leq k\leq n-1$, then all $D_i$'s are Fredholm, and so we trivially choose $A=\mathbf{0}$. Assume the validity of a lengthy condition expressed in $(i)(b)$. Then, one easily checks that one of the Cases 1 or 2 in the proof of \cite[Theorem 3]{WU2} holds, and so we get $A\in\mathcal{B}_n$ so that $T_n^d(A)\in\Phi(\mathcal{H}_1\oplus\cdots\oplus\mathcal{H}_n)$ as described there. $\square$
\begin{Posledica}(\cite[Theorem 3]{WU2}, corrected version)\label{POSLEDICA9}
Let $D_1\in\mathcal{B}(\mathcal{H}_1),\ D_2\in\mathcal{B}(\mathcal{H}_2),...,D_n\in\mathcal{B}(\mathcal{H}_n)$. Then
$$
\begin{aligned}
\sigma_{SF+}(D_1)\cup\sigma_{SF-}(D_n)\cup\Big(\bigcup\limits_{k=2}^{n-1}\Delta_k\Big)\cup\Delta_n\subseteq\\\bigcap_{A\in\mathcal{B}_n}\sigma_e(T_n^d(A))\subseteq\\\sigma_{SF+}(D_1)\cup\sigma_{SF-}(D_n)\cup\Big(\bigcup\limits_{k=2}^{n-1}\Delta_k\Big)\cup\Delta_n\cup\Big(\bigcup_{k=2}^{n-1}\Delta_k'\Big),
\end{aligned}
$$
where
$$
\begin{aligned}
\Delta_k=\lbrace\lambda\in\mathds{C}:\ \alpha(D_k-\lambda)=\infty\ and\ \sum_{s=1}^{k-1}\beta(D_s-\lambda)<\infty\rbrace\cup\\
\lbrace\lambda\in\mathds{C}:\ \beta(D_k-\lambda)=\infty\ and\ \sum_{s=k+1}^n\alpha(D_s-\lambda)<\infty\rbrace,\quad 2\leq k\leq n-1,
\end{aligned}
$$
$$
\begin{aligned}
\Delta_n=\lbrace\lambda\in\mathds{C}:\ \alpha(D_n-\lambda)=\infty\ and\ \sum_{s=1}^{n-1}\beta(D_s-\lambda)<\infty\rbrace\cup\\\lbrace\lambda\in\mathds{C}:\ \beta(D_1-\lambda)=\infty\ and\ \sum_{s=2}^n\alpha(D_s-\lambda)<\infty\rbrace,
\end{aligned}
$$
$$
\Delta_k':=\Big\lbrace\lambda\in\mathds{C}:\ \mathcal{R}(D_k-\lambda)\ is\ not\ closed\Big\rbrace,\ 2\leq k\leq n-1,
$$

\end{Posledica}

\begin{Primedba}
Again, due to the presence of sets $\Delta_k'$, $2\leq k\leq n-1$, we have a difference between definitions of collections $\Delta_k$, $2\leq k\leq n-1$, in Corollary \ref{POSLEDICA9} and in \cite[Theorem 3]{WU2}.  
\end{Primedba}

We get some interesting results for $n=2$ that seem new in the literature.
\begin{Teorema}\label{FREDHOLMn=2}
Let $D_1\in\mathcal{B}(\mathcal{H}_1),\ D_2\in\mathcal{B}(\mathcal{H}_2)$. Consider the following statements:\\
$(i)$ There exists $A\in\mathcal{B}_2$ such that $T_2^d(A)\in\Phi(\mathcal{H}_1\oplus\mathcal{H}_2)$;\\[3mm]
$(ii)$   $(a)$ $D_1\in\Phi_+(\mathcal{H}_1)$ and $D_2\in\Phi_-(\mathcal{H}_2)$;\\
\hspace*{7.5mm}$(b)$ $\beta(D_1)=\alpha(D_2)=\infty$ or \Big($D_2\in\Phi_+(\mathcal{H}_2)$ and $D_1\in\Phi_-(\mathcal{H}_1)$\Big)\\[3mm]
Then $(i) \Leftrightarrow (ii)$.
\end{Teorema}

\begin{Posledica}\label{POSLEDICA10}
Let $D_1\in\mathcal{B}(\mathcal{H}_1),\ D_2\in\mathcal{B}(\mathcal{H}_2)$. Then
$$
\bigcap_{A\in\mathcal{B}_2}\sigma_e(T_2^d(A))=\sigma_{SF+}(D_1)\cup\sigma_{SF-}(D_2)\cup\Delta,
$$
where
$$
\begin{aligned}
\Delta=\lbrace\lambda\in\mathds{C}:\ \alpha(D_2-\lambda)=\infty\ and\ \beta(D_1-\lambda)<\infty\rbrace\cup\\\lbrace\lambda\in\mathds{C}:\ \beta(D_1-\lambda)=\infty\ and\ \alpha(D_2-\lambda)<\infty\rbrace.
\end{aligned}
$$

\end{Posledica}

\noindent\textbf{Acknowledgments}\\[3mm]
\hspace*{6mm}I wish to express my gratitude to professor Dragan S. Djordjevi\'{c} for introduction to this topic and for useful comments that greatly improved the form of this paper.  \\[3mm]

\noindent\textbf{Disclosure statement}\\[3mm]
\hspace*{6mm}I declare there are no conflicts of interest associated with this work.  \\[3mm]

\end{document}